\documentclass{amsart}
\usepackage{graphicx}
\usepackage{amssymb}

\newtheorem{theorem}{Theorem}[section]
\newtheorem{lemma}[theorem]{Lemma}
\newtheorem{corollary}[theorem]{Corollary}

\theoremstyle{definition}
\newtheorem{definition}[theorem]{Definition}
\newtheorem{example}[theorem]{Example}

\theoremstyle{remark}

\newtheorem{question}[theorem]{Question}

\numberwithin{equation}{section}

\begin{document}
\title [Indecomposability from the Complement]{Characterizing
  Indecomposable Plane Continua From Their Complements}

\author[C.~P.~Curry]{Clinton P.~Curry}
\email[Clinton P.~Curry]{clintonc@uab.edu}

\author[J.~C.~Mayer]{John C.~Mayer}
\email[John C.~Mayer]{mayer@math.uab.edu}

\address[Clinton P.~Curry and John C.~Mayer]
{Department of Mathematics\\ University of Alabama at Birmingham\\
  Birmingham, AL 35294-1170}

\author[E.~D.~Tymchatyn]{E.~D.~Tymchatyn}
\email[E.~D.~Tymchatyn]{tymchat@math.usask.ca}

\address[E.~D.~Tymchatyn]{Dept.~of Mathematics and Statistics\\
  University of Saskatchewan\\
  Saskatoon, Saskatchewan, Canada S7N 0W0}

\keywords{indecomposable continuum, complementary domain, Julia set,
complex dynamics, buried point} \subjclass[2000]{Primary: 54F15;
Secondary: 37F20}

\thanks{The third named author was supported in part by NSERC
  0GP005616.  We thank the Department of Mathematics and Computer Science at Nipissing
  University, North Bay, Ontario, for the opportunity to work on this paper in pleasant
  surroundings at their annual topology workshop.}

\date{\today}

\begin{abstract}
  We show that a plane continuum $X$ is indecomposable iff $X$ has a
  sequence $(U_n)_{n=1}^\infty$ of not necessarily distinct complementary domains satisfying
  the {\em double-pass condition}: for any sequence
  $(A_n)_{n=1}^\infty$ of open arcs, with $A_n \subset U_n$ and
  $\overline{A_n}\setminus A_n \subset \partial U_n$, there is a sequence
  of shadows $(S_n)_{n=1}^\infty$, where each $S_n$ is a shadow of
  $A_n$, such that $\lim S_n=X$. Such an open arc divides $U_n$ into
  disjoint subdomains $V_{n,1}$ and $V_{n,2}$, and a \emph{shadow} (of
  $A_n$) is one of the sets $\partial V_{n,i}\cap \partial U$.
\end{abstract}

\maketitle

\section{Introduction}

In this paper, a {\em continuum} is a compact, connected, nonempty
metric space.  A continuum is {\em decomposable} if it can be
written as the union of two of its proper subcontinua; otherwise, it
is {\em indecomposable}.  Let $\mathbb{C}$ denote the complex plane
and let $\mathbb{C}_\infty$ denote the Riemann sphere
$\mathbb{C}\cup\{\infty\}$. A \emph{plane domain} is a subset of
$\mathbb{C}_\infty$ which is conformally isomorphic to the open unit
disk $\mathbb{D} \subset \mathbb{C}_\infty$ (which is to say that it
is open, connected, simply connected and its boundary is a
nondegenerate subcontinuum of $\mathbb{C}_\infty$).  If $X$ is a
continuum in $\mathbb{C}_\infty$, the components of
$\mathbb{C}_\infty\setminus X$ are called \emph{complementary
domains} and are plane domains. If $W\subset \mathbb{C}_\infty$, we
denote the boundary of $W$ by $\partial W$.  We say that a point $x$
of a continuum is \emph{buried} if it does not lie on the boundary
of any complementary domain.  The spherical metric on
$\mathbb{C}_\infty$ is denoted by $d$, and $H_d$ denotes the
Hausdorff metric on the hyperspace of subcontinua of
$\mathbb{C}_\infty$ \cite[Section 4.1]{nadler}.

There are several ways of recognizing intrinsically that a continuum
$X$ is indecomposable. For instance, $X$ is indecomposable if and only
if every proper subcontinuum of $X$ is nowhere dense in $X$
\cite{hockyoun}.  Also, a continuum $X$ is indecomposable if there are
points $a,b,c\in X$ such that no proper subcontinuum of $X$ contains
any two of these points.  In this paper, we are interested in
recognizing indecomposable planar continua not by intrinsic
properties, but by the relationship between the continuum and its
ambient space.

Indecomposable continua arise naturally in dynamical systems
\cite{deva04,devajarq02}.  However, in specific dynamical systems,
it is often difficult to recognize them.  In complex analytic
dynamics, the {\em Julia set} of a rational map
$f:\mathbb{C}_\infty\to\mathbb{C}_\infty$ is the set of unstable
points under iteration of $f$ (see \cite{miln} for definitions). A
long-standing question in complex analytic dynamics asks: {\em Can
the Julia set of a rational function be an
  indecomposable continuum?} Several authors have attacked this
question, among them \cite{mayeroge93, ChMaRo03} for polynomials,
\cite{yangsun03} for {\em bicritical} rational maps (rational maps
with exactly two critical points), and \cite{ChMaTYTu05} for Julia
sets of a class of rational functions with no buried points.  In this
situation, it is much easier to analyze the complement of the Julia
set, called the \emph{Fatou set}; this motivates our interest in
studying indecomposability from the point of view of a continuum's
complement.

The second named author, with various co-authors
\cite{mayeroge93,ChMaRo03,ChMaTYTu05}, investigated the recognition
of indecomposable continua from their complement in the case that
$\partial U = X$ for some complementary domain $U$ of $X$. The tool
used was prime end theory, and a characterization was obtained in
the context of $X$ being a Julia set, making use of the dynamics.
The characterization of indecomposable continua from their
complements in the current paper primarily addresses the case that
$X$ is not the boundary of any of its complementary domains (which
would imply that there are infinitely many complementary domains,
each having boundaries nowhere dense in $X$). Even better, this
characterization also subsumes the first case and is entirely
topological.

To state our characterization theorem, we need some definitions.
These concepts are related to those which arose originally in prime
end theory.

\begin{definition}
  Let $U$ be a plane domain.  A {\em crosscut} of $U$ is an open arc
  $A=(a,b)\subset U$ such that $\overline A = [a,b]$ is a closed arc which
  meets $\partial U$ exactly in the set $\{a,b\}$.  A \emph{generalized
    crosscut} of $U$ is an open arc $A\subset U$ such that $\overline
  A \setminus A \subset \partial U$.
\end{definition}

Notice that the notion of a generalized crosscut is strictly broader
than the notion of a crosscut. It is easy to see that a generalized
crosscut of a domain $U$ cuts $U$ into two nonempty disjoint
subdomains $V_1$ and $V_2$ such that $U=V_1\cup A\cup V_2$.

\begin{definition}
  Let $U$ be a plane domain and $A$ a generalized crosscut of $U$.  We
  call each component of $U \setminus A$ a {\em crosscut
    neighborhood}.  If $V$ is a crosscut neighborhood determined by
  a generalized crosscut $A$, we call the continuum $S=\partial V\cap \partial U$
  a {\em shadow} of $A$.
\end{definition}

Thus, a generalized crosscut $A$ of a domain $U$ has exactly two
crosscut neighborhoods, and consequently two shadows whose union is
$\partial U$.  Examples below show that one or both of these shadows
can be proper subcontinua of $\partial U$ or, more surprisingly, all
of $\partial U$.

Limits below are interpreted in the metric $H_d$.

\begin{definition}\label{double-pass}
  A sequence $(U_n)_{n=1}^\infty$ of (not necessarily distinct)
  complementary domains of a continuum $X$ satisfies the
  \emph{double-pass condition} if, for any sequence of generalized
  crosscuts $A_n$ of $U_n$, there is a sequence of shadows
  $(S_n)_{n=1}^\infty$ of $(A_n)_{n=1}^\infty$ such that $\lim_{n
    \rightarrow \infty} S_n = X$.
\end{definition}

In Section \ref{sec:char}, we prove the following theorem, which is
the main theorem of this paper.

\begin{theorem}[Characterization Theorem] \label{char}
  A planar continuum $X$ is indecomposable if and only if it has a
  sequence $(U_n)_{n=1}^\infty$ of complementary domains which
  satisfies the double-pass condition.
\end{theorem}

\section{Partial and Prior Results}
\subsection{Brief History}
The first partial recognition theorem for indecomposable continua from
the complement is that of Kuratowski \cite{kura}.

\begin{theorem}[Kuratowski]\label{theorem:Kuratowski}
  If a plane continuum $X$ is the common boundary of three of its
  complementary domains, then $X$ is either indecomposable or the
  union of two proper indecomposable subcontinua.
\end{theorem}

The following theorem of Rutt seems quite different, and is also only
applicable if $X$ is the boundary of some complementary domain.

\begin{theorem}[Rutt, \cite{rutt}]\label{theorem:Rutt}
  If a nondegenerate plane continuum $X$ is the boundary of a
  complementary domain $U$, and if there is a prime end of $U$ whose
  impression is $\partial U=X$, then $X$ is either indecomposable or the
  union of two proper indecomposable subcontinua.
\end{theorem}

Without going into detail (but see \cite{ChMaRo03}), the {\em
impression} of a prime end of $U$ is the intersection of the shadows
of a sequence $(A_n)_{n=1}^\infty$ of crosscuts of $U$ having the
property that for each $n$, $(A_m)_{m>n}$ is a pairwise closure
disjoint null sequence contained in one of the crosscut neighborhoods
of $A_n$.

The connection among the theorems above is made explicit by a
technical theorem of Burgess.  While the original result is stated
in terms of what Burgess calls simple disks, the theorem can be
equivalently stated in terms of  closed balls.   For $a \in
\mathbb{C}_\infty$ and $r > 0$, define the {\em ball} of radius $r$
about $a$ by
\[
B_r(a) = \{z \in \mathbb{C}_\infty \mid d(a, z) < r\}.
\]

\begin{theorem}[Burgess, {\cite[Theorem 9]{burg}}]\label{threeD}
  Let $H$ be a closed set and $X$ a continuum in the plane.  Suppose
  $X_1$, $X_2$, and $X_3$ are subcontinua of $X$ and $D_1$, $D_2$, and
  $D_3$ are pairwise disjoint closed balls with $D_i \cap H =
  \emptyset$ and $\emptyset \neq D_i \cap X_i = D_i \cap X$ for each
  $i \in \{1,2,3\}$. Then there do not exist three distinct
  complementary domains of $X\cup H$ such that each of them intersects
  each of the balls $D_i$.
\end{theorem}

Using this theorem, Burgess proves the following recognition theorem,
which also applies when the continuum is not the union of the
boundaries of its complementary domains.

\begin{corollary}[Burgess, {\cite[Corollary to Theorem 9]{burg}}]\label{Burgess}
  If the plane continuum $X$ is the limit of a sequence of distinct
  complementary domains of $X$, then either $X$ is indecomposable, or
  there is only one pair of indecomposable continua whose union is
  $X$.
\end{corollary}

As recognition theorems, the above suffer from the weakness of their
conclusion.  In \cite{ChMaRo03,ChMaTYTu05} dynamical considerations
rule out that the Julia set of a polynomial can be the union of two
proper indecomposable subcontinua.  However, this is under the
hypothesis that the Julia set is the boundary of one of its
complementary domains.  The following definition and recognition
theorem appear in \cite{ChMaTYTu05}.  Since it represents a
simplification of the proof in \cite{ChMaTYTu05}, we prove
Theorem~\ref{theorem:antichain} making use of Theorem~\ref{threeD}.

\begin{definition}
  An \emph{antichain} of crosscuts of a plane domain $U$ is a sequence
  $(H_n)_{n=1}^\infty$ of distinct pairwise closure disjoint crosscuts
  of $U$ such that, for each $m$, one crosscut neighborhood of $H_m$
  contains all the crosscuts $(H_n)_{n \neq m}$.
\end{definition}

\begin{theorem}\label{theorem:antichain}
  Let $U \subset \mathbb{C}_\infty$ be a plane domain.  Let $z \in U$.  Suppose
  there exists an antichain $(H_n)_{n=1}^\infty$ of crosscuts of $U$
  such that $\lim_{n \rightarrow \infty} \mathrm{Sh}(H_n)=\partial U$, where
  $\mathrm{Sh}(H_n)$ is the shadow of the crosscut neighborhood $W_n$ of $U
  \setminus H_n$ which misses $z$.  Then $\partial U$ is indecomposable.
\end{theorem}

\begin{proof}
  For a contradiction, suppose $\partial U$ satisfies the hypotheses of the
  theorem, but may be written as the union of proper subcontinua $X_1$
  and $X_2$.  By passing to a subsequence, we may assume that
  $(H_n)_{n=1}^\infty$ converges to a point of $\partial U$.  Choose disjoint
  closed balls $D_1$ and $D_2$ such that
  \begin{enumerate}
    \item $z \notin D_i$,
        \item $D_i \cap H_n = \emptyset$ for all $i \in \{1,2\}, n
          \in \mathbb{N}$,
    \item $X_i$ intersects the interior of $D_i$, and
    \item $D_i \cap X_j = \emptyset$ if $i \neq j$.
  \end{enumerate}
  Choose three crosscuts $H_1$, $H_2$, and $H_3$ so that the component
  $W_i$ of $\mathbb{C}_\infty \setminus (\partial U \cup H_i)$ missing $z$ hits both
  $D_1$ and $D_2$. Notice that, since the crosscuts are members of an
  antichain, $W_i \cap W_j =\emptyset$ for distinct $i$ and $j$ in
  $\{1,2,3\}$. Let $R_1$, $R_2$, and $R_3$ be arcs from $z$
  to $\partial U$ disjoint (except for
  $z$) from each other and from $D_1\cup D_2$, and such that each $H_i$ lies in a different complementary
  component $U_i$ in $U$ of
    \[ X = \partial U \cup R_1 \cup R_2 \cup R_3. \]
  Notice that $W_i \subset U_i$. Define $X_3 = R_1 \cup R_2 \cup R_3$,
  and let $D_3 \subset U$ be a closed ball about $z$ which is disjoint
  from $D_1$ and $D_2$. Thus, each $U_i$ intersects each of $D_1$, $D_2$, and
  $D_3$; this contradicts Theorem \ref{threeD}, with $H=\emptyset$
 in the statement.
\end{proof}

\subsection{Necessary Condition}\label{necessity}

In this section we show that for a plane continuum $X$ to be
indecomposable, it is necessary that $X$ have a sequence of
complementary domains whose boundaries converge to $X$.  The proof
requires a few additional facts and definitions.

\begin{definition}
  The {\em composant}, denoted $C(p)$, of a point $p$ in a continuum
  $X$ is the union of all the proper subcontinua of $X$ that contain
  $p$.
\end{definition}

\begin{theorem}[{\cite{hockyoun}}] \label{theorem:hockyoun}
  Let $X$ be a nondegenerate indecomposable continuum.  Then the
  following hold:
  \begin{enumerate}
    \item $X$ has $c$ pairwise disjoint composants.
    \item Each composant is dense in $X$.
    \item Each composant can be written as a countable increasing
      union of nowhere dense proper subcontinua of $X$, converging to
      $X$ in the Hausdorff metric.
\end{enumerate}
\end{theorem}

\begin{definition}
  A connected topological space $X$ is said to be \emph{unicoherent}
  if, for any pair $A$ and $B$ of closed, connected subsets such that
  $A \cup B = X$, the intersection $A \cap B$ is connected.
\end{definition}

Note that the plane itself and an open or closed ball in the plane
is unicoherent \cite{stone1949}.  Recall that $B_r(a)$ denotes the
open ball of radius $r>0$ about center $a$.

\begin{theorem}\label{necplane}
  Let $X$ be an indecomposable plane continuum.  Then there is a
  sequence $(U_n)_{n=1}^\infty$ of (not necessarily distinct)
  complementary domains of $X$ such that $\lim\partial U_n=X$.
\end{theorem}

\begin{proof}
  This is clear if $X$ is a point, so assume $X$ is a nondegenerate
  indecomposable continuum.  Take $p,q,r \in X$, each in a different
  composant of $X$.  For each $n\in\mathbb{N}$, define
  \begin{align*}
    Q_n &= \text{the component of } X \setminus B_{1/n}(p) \text{
      containing } q\\
    R_n &= \text{the component of } X \setminus B_{1/n}(p) \text{
      containing } r
  \end{align*}

  Notice that $\lim_{n\rightarrow \infty}Q_n = \lim_{n \rightarrow
  \infty} R_n = X$, by Theorem \ref{theorem:hockyoun}.  Since $Q_n$ and
  $R_n$ are different components of $X
  \setminus B_{1/n}(p)$, they are separated in $\mathbb{C}_\infty \setminus
  B_{1/n}(p)$ by $\mathbb{C}_\infty \setminus
  (B_{1/n}(p) \cup X)$.  Also, $Q_n$ and $R_n$ are closed in the
  normal space $\mathbb{C}_\infty \setminus B_{1/n}(p)$, so there is a subset
  $K_n$, closed in $\mathbb{C}_\infty \setminus B_{1/n}(p)$, of
  $\mathbb{C}_\infty \setminus (B_{1/n}(p) \cup X)$ which separates $Q_n$ and
  $R_n$.  Since $\mathbb{C}_\infty \setminus B_{1/n}(p)$ is homeomorphic to the
  closed unit disk in the plane, it is unicoherent, so a component
  $L_n$ of $K_n$ is a closed (in $\mathbb{C}_\infty \setminus B_{1/n}(p)$) and
  connected subset of $\mathbb{C}_\infty \setminus (B_{1/n}(p)\cup X)$ which
  separates $Q_n$ and $R_n$ in $\mathbb{C}_\infty \setminus B_{1/n}(p)$
  \cite{stone1949}. Moreover, since $L_n \subset \mathbb{C}_\infty \setminus
  X$, it lies in a single complementary domain $U_n$ of $X$. The
  sequence $(U_n)_{n=1}^\infty$ formed in this way is the required
  sequence of complementary domains.

  It is evident that $\lim_{n \rightarrow \infty} \partial U_n \subset X$;
  we aim to show that $X \subset \lim_{n \rightarrow \infty} \partial U_n$.
  Choose $\epsilon>0$, and $x \in X$.  Let $N \in \mathbb{N}$ such that, for
  all $n \ge N$
  \begin{enumerate}
    \item $Q_n \cap B_\epsilon(x) \neq \emptyset$,
    \item $R_n \cap B_\epsilon(x) \neq \emptyset$, and
    \item $B_{1/n}(p) \cap B_\epsilon(x) = \emptyset$.
  \end{enumerate}

  For $n \ge N$, choose $q_n \in Q_n \cap B_\epsilon(x)$ and $r_n \in
  R_n \cap B_\epsilon(x)$ for $n \ge N$. The straight line segment
  $A_{n}$ from $q_n$ to $r_n$ is a subset of $B_\epsilon(x)$ and,
  hence, of $\mathbb{C}_\infty \setminus B_{1/n}(p)$, so $A_{n}$ must meet $L_n$,
  since $L_n$ separates $q_n$ from $r_n$ in $\mathbb{C}_\infty \setminus
  B_{1/n}(p)$. Since $L_n\subset U_n$, $A_{n}$ intersects $U_n \cap
  B_\epsilon(x)$. Since $q_n$, $r_n$ are not in $U_n$ (they lie in
  $X$), $A_{n}$ intersects $\partial U_n$, and $\partial U_n \cap
  B_{\epsilon}(x) \neq \emptyset$. This is true for all $n \ge N$, so
  $x \in \liminf_{n \rightarrow \infty} \partial{U_n} \subset \lim_{n
  \rightarrow \infty}\partial{U_n}$. This completes the proof.
\end{proof}

\section{The Characterization Theorem} \label{sec:char}
We saw in Subsection~\ref{necessity} that having a sequence of
complementary domains whose boundaries converge to $X$ is a necessary
condition for the plane continuum $X$ to be indecomposable.
Example~\ref{dblK} below shows that this condition is not sufficient,
even if the domains are distinct, and suggests that we must find a way
to rule out that the sequence of complementary domains ``splits'' into
``halves'' each of which converge to proper indecomposable subcontinua.

\subsection{Examples}

\begin{figure}
  \includegraphics[width=2in]{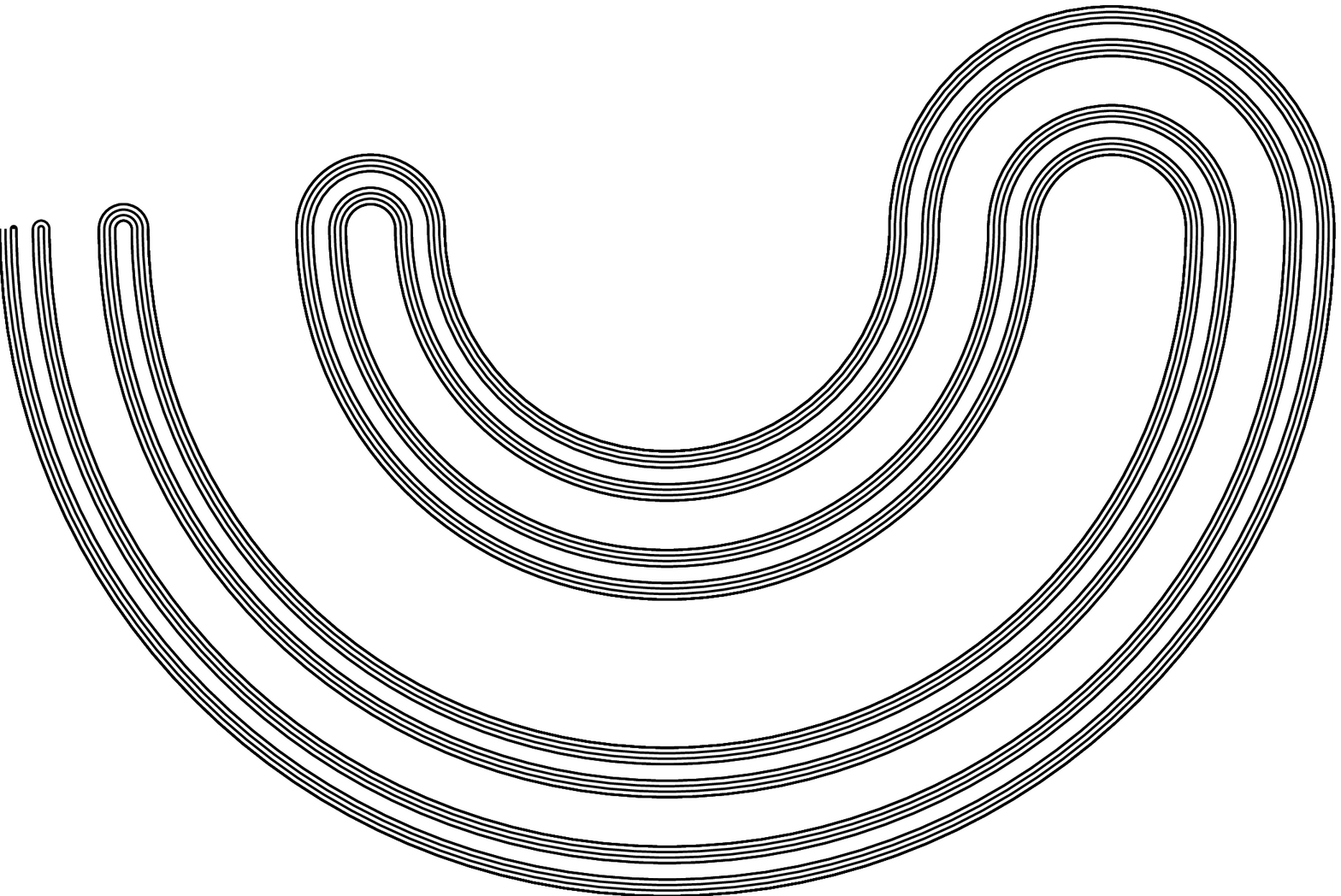}
  \includegraphics[width=2in]{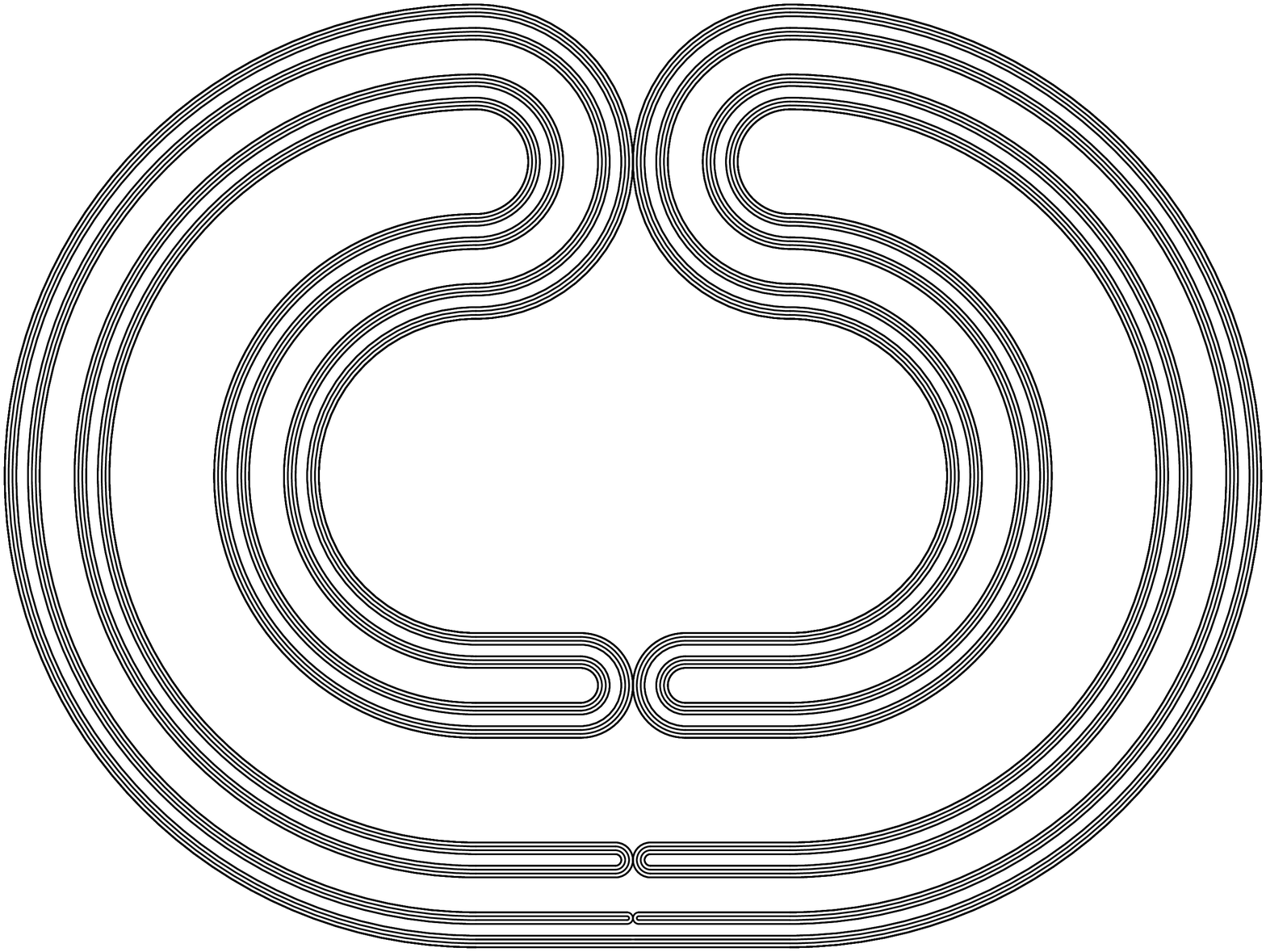}
  \caption{The Knaster buckethandle continuum (left); a union of two
    Knaster continua (right) meeting at a sequence of points converging to their common endpoint.}
  \label{snglK_fig}
\end{figure}

\begin{example}\label{snglK}
  The Knaster buckethandle continuum, depicted on the left in
  Figure~\ref{snglK_fig}, is a standard example of an indecomposable
  continuum.  It can be viewed as a disk from which successively deeper
  and denser fjords are dug.  Notice that, for any generalized crosscut
  drawn in its single complementary domain, infinitely many fjords lie
  in one crosscut neighborhood or the other, so one shadow is dense.
  Notice that one composant of $X$ is the one-to-one continuous image of
  the half line $[0,\infty)$.  The point of $X$ corresponding to the
  point $0$ of $[0, \infty)$ is called the endpoint of $X$.  For a more
  precise construction, see \cite[Vol. II, p. 205]{kura}.
\end{example}

\begin{example}\label{dblK}
  The continuum $X$ depicted on the right in Figure~\ref{snglK_fig} is
  an example of a continuum which satisfies the hypotheses of
  Theorem~\ref{Burgess} without being indecomposable.  It is a
  symmetric union of two Knaster buckethandle continua $X_1$ and $X_2$
  intersecting in a countable set which lies on a vertical line.  The
  continuum has infinitely many complementary domains
  $(U_i)_{i=1}^\infty$.  The decomposability of $X$ can be detected
  from the complementary domains as follows: to its complementary
  domains $(U_i)_{i=1}^\infty$, associate the collection of crosscuts
  $(K_i)_{i=1}^\infty$, where $K_i$ lies in $U_i$ on the vertical axis
  of symmetry.  Each crosscut has one shadow which is a subcontinuum
  of $X_1$, and another which is a subcontinuum of $X_2$.  Therefore,
  any convergent sequence of shadows must limit to a proper
  subcontinuum of $X$.
\end{example}

\begin{figure}
  \includegraphics[width=2in]{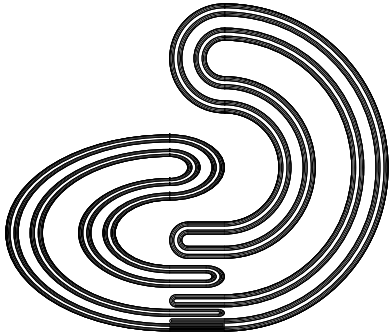}
  \caption{Two Knaster continua meeting on an end interval of each.}
  \label{nocrosscuts_fig}
\end{figure}

\begin{example}\label{nocrosscuts}
  The continuum in Figure~\ref{nocrosscuts_fig} is the union of a pair
  of Knaster continua $X_1$ and $X_2$ with distinct endpoints such that
  $X_1 \cap X_2$ is the horizontal arc $A$ between the endpoints of
  $X_1$ and $X_2$.  This continuum, like Example \ref{snglK}, has the
  property that every crosscut has a dense shadow, despite the
  continuum's decomposability. Let $K$ be a generalized crosscut such
  that one end lands on a point of $X \setminus A$ and the other end
  wiggles between the Knaster continua and compactifies on $A$.  Neither
  shadow of this generalized crosscut is dense, so the constant sequence
  consisting of this crosscut fails the double pass condition.
\end{example}

\subsection{Proof of Characterization Theorem}

In Definition \ref{double-pass}, we defined the double-pass condition
on a sequence of generalized crosscuts in a sequence of complementary
domains which is motivated by Example~\ref{dblK} and by the
similarly-functioning condition of Cook and Ingram \cite{cookingr}
introduced for recognizing indecomposable (chainable) continua in
terms of refining open covers.  Here we prove our main theorem: the
existence of a sequence of complementary domains satisfying our
double-pass condition is equivalent to indecomposability.

The following Lemma follows from \cite[(A1.4)]{ursellyoung}, but we
include a self-contained proof here for convenience.

\begin{lemma}\label{lem:nullsequencecrosscuts}
  Suppose that $\phi:U \rightarrow \mathbb{D}$ is a conformal isomorphism,
  where $U$ is a plane domain.  Then the image of a null sequence
  $(K_n)_{n=1}^\infty$ of crosscuts of $U$ is a null sequence of
  crosscuts of $\mathbb{D}$.
\end{lemma}
\begin{proof}
  By way of contradiction, let $(K_n)_{n=1}^\infty$ be a null sequence
  of crosscuts of $U$ such that the image sequence $(A_n)_{n=1}^\infty
  = (\phi(K_n))_{n=1}^\infty$ consists of crosscuts whose diameters
  are bounded away from zero. Then, by passing to a subsequence, we
  may assume that the image sequence accumulates on a nondegenerate
  continuum $L \subset \overline{\mathbb{D}}$.  Since $(K_n)_{n=1}^\infty$
  does not accumulate on a subset of $U$, we see that $L \subset \partial
  \mathbb{D}$.  Also, we may assume without loss of generality that
  $(K_n)_{n=1}^\infty$ converges to a point of $x \in \partial U$.

  Let $t \in L$. There exists a chain of crosscuts
  $(A_n')_{n=1}^\infty$ of $\mathbb{D}$ converging to $t$ which maps to a
  null sequence $(K_n')_{n=1}^\infty$ of crosscuts of $U$ by
  $\phi^{-1}$ (see \cite[Lemma 17.9]{miln}).  We may assume that
  $(K_n')_{n=1}^\infty$ converges to a point of $\partial U$ by passing to
  a subsequence.  Since $(A_n)_{n=1}^\infty$ accumulate on $t$, all
  but finitely many $A_n$ intersect the crosscut neighborhood of
  $A_m'$ corresponding to $t$.  Also, since $(K_n)_{n=1}^\infty$ forms
  a null sequence in $U$, we see that all but finitely many $K_n$
  (thus $A_n$) lie entirely within the crosscut neighborhood of $K_m'$
  (thus $A_m'$) corresponding to $t$.  However, the crosscut
  neighborhoods of $A_m'$ form a null sequence, so
  $(A_n)_{n=1}^\infty$ form a null sequence.
\end{proof}

We say that a pair of subsets $E_1$ and $E_2$ of
$\partial\mathbb{D}$ are $\emph{unlinked}$ if there exist intervals
$I_1 \supset E_1$ and $I_2 \supset E_2$ such that $|I_1 \cap I_2|
\le 2$.

\begin{lemma}\label{lem:unlinked}
  Suppose $U \subset \mathbb{C}_\infty$ is a plane domain, and let $\phi:U
  \rightarrow \mathbb{D}$ be a conformal isomorphism.  Let $B_1$ and $B_2$
  be disjoint closed  balls meeting $\partial U$ in their interiors, and let $E_i
  \subset \partial\mathbb{D}$ denote the endpoints of all of the crosscuts of $\mathbb{D}$
  which constitute $\phi( (\partial B_i) \cap U)$.  If $E_1$ and $E_2$ are
  unlinked, then there is a generalized crosscut $K$ of $U$ which
  separates $B_1 \cap U$ from $B_2 \cap U$ in $U$.  Moreover, if $\partial
  U$ is locally connected, $K$ is a crosscut of $U$.
\end{lemma}

\begin{proof}
  Let $I_1$ and $I_2$ be minimal closed intervals of $\partial\mathbb{D}$ such
  that $I_i \supset E_i$ for $i \in \{1,2\}$ and $|I_1 \cap I_2| \le
  2$.  Note that $(\partial B_i) \cap U$ is the union of a null sequence of
  crosscuts of $U$, so $\phi( (\partial B_i) \cap U)$ is the union of a
  null sequence of crosscuts of $\mathbb{D}$ by
  Lemma~\ref{lem:nullsequencecrosscuts}.  $I_1$ and $I_2$ are unlinked
  implies that neither $\phi(B_1 \cap U)$ nor $\phi(B_2 \cap U)$
  separates the other in $\mathbb{D}$.  By connectedness of $\mathbb{D}$, there
  is then a unique component $D$ of $\mathbb{D} \setminus \phi((B_1 \cup
  B_2) \cap U)$ which meets both $\phi(\partial B_1 \cap U)$ and $\phi(\partial
  B_2 \cap U)$. The crosscuts on $\partial D$ form a null sequence, so it
  is not difficult to show that $\overline D \subset \overline{\mathbb{D}}$
  is locally connected.

  The endpoints of $I_1$ and $I_2$ are each on $\partial D$. Let $K$ be a
  crosscut joining an endpoint of $I_1$ to an endpoint of $I_2$ such
  that $K$ separates the interiors of $I_1$ and $I_2$, in $\partial\mathbb{D}$,
  and thus $\phi(B_1)$ and $\phi(B_2)$ in $\mathbb{D}$.  Then
  $\phi^{-1}(K)$ is a generalized crosscut of $U$ which separates $B_1
  \cap U$ from $B_2 \cap U$.  Further, if $\partial U$ is locally
  connected, $\phi^{-1}$ extends to a continuous function $\overline
  \phi^{-1} : \overline{\mathbb{D}} \rightarrow \overline U$ \cite[Theorem
  17.14]{miln}. In this case, $\phi^{-1}(K)$ is a true crosscut of
  $U$, as $\overline{\phi^{-1}(K)} = \overline{\phi}^{-1}(\overline
  K)$ is an arc.
\end{proof}

Now we have the tools to prove our Characterization
Theorem~\ref{char}.

\begin{proof}[Proof of Theorem~\ref{char}]
  First, suppose that $X$ is indecomposable.  We show $X$ satisfies
  the double-pass condition.  By Theorem \ref{necplane}, there exists
  a sequence $(U_n)_{n=1}^\infty$ of complementary domains of $X$ such
  that $\lim_{n\rightarrow \infty} \partial U_n=X$.  Let
  $(K_n)_{n=1}^\infty$ be a sequence of generalized crosscuts, with
  $K_n$ in $U_n$ for each $n \in \mathbb{N}$.  Let $A_n$ and $B_n$ be the
  shadows of $K_n$, with $H_d(A_n,X) \le H_d(B_n,X)$, where $H_d$
  denotes the Hausdorff metric.

  We claim that $\lim_{n\rightarrow \infty} A_n = X$.  Since the
  hyperspace of subcontinua of $X$ is a compact metric space, it is
  sufficient to show that every convergent subsequence of
  $(A_n)_{n=1}^\infty$ converges to $X$.  Let $(n_i)_{i=1}^\infty$ be
  such that $(A_{n_i})_{i=1}^\infty$ converges to a continuum $A
  \subset X$. By passing to a subsequence, we may assume
  $(B_{n_i})_{i=1}^\infty$ also converges to a continuum $B \subset
  X$.  Since $\lim_{i \rightarrow \infty} \partial U_{n_i} = X$ and
  $A_{n_i} \cup B_{n_i} = \partial U_{n_i}$, we have that $A \cup B = X$.
  Since $X$ is indecomposable, not both $A$ and $B$ may be proper
  subcontinua of $X$, so $A=X$ or $B=X$.  Since,
  $H_d(A_{n_i},X) \le H_d(B_{n_i},X)$ for all $i$, we have
  $A=X$.  This concludes the proof of this implication.

  Now we prove the converse. Let $X$ be a continuum with a sequence
  $(U_n)_{n=1}^\infty$ of complementary domains satisfying the
  double-pass condition.  Suppose, by way of contradiction, that $X =
  X_1 \cup X_2$, where $X_1$ and $X_2$ are proper subcontinua of $X$.
  We can then find  open balls $B_1$ and $B_2$ such that
  \begin{enumerate}
    \item $\overline B_1 \cap \overline B_2 = \emptyset$,
    \item $ B_i \cap X_i\not=\emptyset$ for $i=1,2$,
      and
    \item $\overline B_i \cap X_j = \emptyset$ for $i\not=j$.
  \end{enumerate}

  Since $(U_n)_{n=1}^\infty$ satisfies the double-pass condition,
  there exists a particular $N \in \mathbb{N}$ such that, for any
  generalized crosscut $K$ of $U_N$, one shadow of $K$ intersects both
  $B_1$ and $B_2$.  We fix $U = U_N$, and let $\phi:U \rightarrow
  \mathbb{D}$ be a conformal isomorphism.  Define $E_1$ and $E_2$, as in
  Lemma~\ref{lem:unlinked}, to be the sets of endpoints of the crosscuts
  comprising $\phi( (\partial B_1) \cap U)$ and $\phi( (\partial B_2) \cap U)$,
  respectively.  There are two cases: Either $E_1$ and $E_2$ are
  linked, or they are not.  The second case cannot occur, as
  Lemma~\ref{lem:unlinked} asserts the existence of a generalized
  crosscut $K_0$ of $U$ separating $B_1 \cap U$ and $B_2 \cap U$,
  contrary to our assumption.

  Thus, $E_1$ and $E_2$ are linked.  Note that each of
  $\phi( (\partial B_1 )\cap U)$ and $\phi( (\partial B_2 )\cap U)$ consists
  of crosscuts of $\mathbb{D}$ with endpoints in $E_1$ and $E_2$,
  respectively.
  There are two cases: (1) either
  one of $\phi( (\partial B_1 )\cap U)$ or $\phi( (\partial B_2 )\cap U)$
  separates the other in $\mathbb{D}$, or (2) neither $\phi( (\partial B_1) \cap U)$ nor
  $\phi( (\partial B_2) \cap U)$ separates the other in $\mathbb{D}$.
 In each case, we  construct a crosscut
  $A\subset \mathbb{D}$. This crosscut will have the
  property that $\phi^{-1}(A)$ is a crosscut of $U$ which
  we show leads to a separation of one of $X_1$ or $X_2$, a contradiction.

  In case (1), without loss of generality, $\phi((\partial B_1) \cap U)$ separates
  $\phi((\partial B_2) \cap U)$ in $\mathbb{D}$.
  Since $\mathbb{D}$ is unicoherent, a component
  of $\phi((\partial B_1) \cap U)$ also separates $\phi((\partial B_2) \cap U)$, so a
  crosscut in $\phi((\partial B_1) \cap U)$ does.  Let $A$ be this
  crosscut.  Then $\phi^{-1}(A)$ is a crosscut of $U$ separating
  $B_2 \cap U$ in $U$.

  For case (2), we suppose that neither $\phi((\partial B_1) \cap U)$ nor
  $\phi((\partial B_2) \cap U)$ separates the other in $\mathbb{D}$.  Since
  $E_1$ and $E_2$ are linked, let $e_1$ and $e_1'$ be points of
  $E_1$ separated in $\partial\mathbb{D}$ by points of $E_2$.  Let $K_1$ and
  $K_1'$ be crosscuts of $\mathbb{D}$ in $\phi((\partial B_1) \cap U)$ that
  have $e_1$ and $e_1'$, respectively, as endpoints.  Since $\phi((\partial B_2) \cap U)$
   does not separate $K_1$ from $K_1'$ in $\mathbb{D}$, there is an arc
   $C$ from a point of $K_1$ to a point of $K_1'$ in
   $\mathbb{D}\setminus\phi(\overline B_2\cap U)$. Let $A \subset K_1 \cup K_1' \cup
  C$ be a crosscut of $\mathbb{D}$ from $e_1$ to $e_1'$ which then separates $\phi(B_2 \cap U)$ in $\mathbb{D}$.
  Because $\phi^{-1}(K_1)$ and $\phi^{-1}(K_1')$ are crosscuts of $U$, we see
  that $\phi^{-1}(A)$ is a crosscut of $U$.   Moreover, $\phi^{-1}(A)$
  separates $B_2 \cap U$ in $U$.

The proof in cases (1) and (2) now proceeds together.
  Let $S_1$ be an irreducible arc which joins points of $\phi((\partial B_2)
  \cap U)$ which are separated by $A$; we may stipulate that $S_1$
  intersects $A$ exactly once, transversely. By applying $\phi^{-1}$
  to both $A$ and $S_1$, we obtain a crosscut $A'$ of $U$ and a
  compact arc $S_1' \subset U$ between points of $\partial B_2$ which
  intersects $A'$ once transversely.  Let $S_2'$ be a compact arc in
  $\overline B_2$ which joins the endpoints of $S_1'$.  Then $S_1' \cup S_2' =
  S$ is a simple closed curve.  Observe that $S \cap X_1 = \emptyset$,
  since $S_1' \subset U$ and $S_2' \subset \overline B_2$.  However, the compact
  arc $\overline{A'}$ joins points of $X_1$ and intersects $S$ exactly
  once, transversely.  Thus, some point of $X_1$ lies inside and
  another point outside of $S$, while $X_1 \cap S = \emptyset$,
  contradicting the connectedness of $X_1$.
\end{proof}

Examination of the proof of Theorem~\ref{char} gives a stronger
theorem for continua whose complementary domains have locally
connected boundaries.

\begin{definition}
  A sequence of complementary domains $(U_n)_{n=1}^\infty$ satisfies
  the \emph{crosscut condition} if, for every sequence of crosscuts
  $(A_n)_{n=1}^\infty$, $A_n \subset U_n$, there exists a choice of
  shadows $S_n$ of $A_n$ such that $\lim S_n = X$.
\end{definition}

Example~\ref{nocrosscuts} showed that this is a strictly weaker
condition than the double-pass condition.  However, the following
shows that, for a certain class of continua, the two notions are
equivalent.

\begin{corollary}
  A planar continuum $X$ whose complementary domains have locally
  connected boundaries is indecomposable if and only if it has a
  sequence $(U_n)_{n=1}^\infty$ of complementary domains which
  satisfies the crosscut condition.
\end{corollary}

This follows  from the proof of Theorem \ref{char}, since the
generalized crosscut of $U$ constructed in the proof with
Lemma~\ref{lem:unlinked} is a crosscut if $\partial U$ is locally
connected.

\section{Questions and Further Results}

We close with a question about rational Julia sets for which our
Characterization Theorem may prove useful, and two theorems by the
first author that will appear in a subsequent paper extending our
results to surfaces.

\begin{question}
  Let $J=J(R)$ be the Julia set of a rational function
  $R:\mathbb{C}_\infty\to\mathbb{C}_\infty$ and suppose that $J$ has buried points.  Can
  $J$ be the union of two proper indecomposable subcontinua?  In
  particular, can $J$ contain a proper indecomposable subcontinuum
  with interior in $J$?
\end{question}

\begin{definition}
  A {\em surface} is a  connected Hausdorff space with a countable
  basis each point of which has a neighborhood homeomorphic to an open
  ball in the plane.  Let $X$ be a continuum in the surface $S$.  As
  before, a component of $S\setminus X$ is called a {\em complementary
  domain}.
\end{definition}

\begin{definition}
  A connected topological space $X$ is \emph{multicoherent} of degree
  $k$ if, for any pair of closed, connected sets $A$ and $B$ such that
  $A \cup B = X$, the intersection $A \cap B$ consists of at most $k$
  components.
\end{definition}

A complementary domain in a surface, unlike in the planar case, need
not be simply connected.  Using the notion of multicoherence and its
consequences (see \cite[Theorem 1]{stonemult1949} for the relevant
extension of the Phragm\`{e}n-Brouwer theorem), we can prove the
following theorem.  We omit the proof, which is similar to the proof
of Theorem~\ref{necplane}.

\begin{theorem}\label{surface}
  Let $S$ be a compact surface and $X$ an indecomposable subcontinuum
  of $S$.  Then there is a sequence $(U_n)_{n=1}^\infty$ of
  complementary domains of $X$ such that $\lim\partial U_n=X$.
\end{theorem}

We claim in Theorem~\ref{surface} that having a sequence of
complementary domains converging to $X$ is a necessary condition for
continuum $X$ contained in a surface $S$ to be indecomposable.  We
saw in the plane a partial converse: given a sequence of distinct
complementary domains $(U_n)_{n=1}^\infty$ such that $\lim\partial
U_n=X$, it follows that $X$ is either indecomposable or the union of
two proper indecomposable subcontinua (Theorem~\ref{Burgess}). In
this connection, we close with the following two theorems
generalizing Burgess's Theorem~\ref{Burgess} and our
Characterization Theorem~\ref{char} to continua in surfaces, proofs
of which will appear subsequently in a paper by the first-named
author.

\begin{theorem}
  Let $S$ be a compact surface.  Suppose continuum $X\subset S$ has a
  sequence $(U_n)_{n=1}^\infty$ of distinct complementary domains with
  $\lim\partial U_n=X$. Then either $X$ is indecomposable, or there is only
  one pair of indecomposable  subcontinua  whose union is $X$.
\end{theorem}

\begin{theorem}
  Let $S$ be a compact surface.  Suppose continuum $X\subset S$ has
  buried points.  Then $X$ is indecomposable iff $X$ has a sequence
  $(U_n)_{n=1}^\infty$ of distinct complementary domains satisfying
  the {\em double-pass condition}: for any sequence
  $(A_n)_{n=1}^\infty$ of generalized crosscuts (suitably defined),
  with $A_n \subset U_n$, there is a sequence of shadows
  $(S_n)_{n=1}^\infty$, where each $S_n$ is a shadow of $A_n$, such
  that $\lim S_n=X$.
\end{theorem}

\bibliographystyle{plain}

\end{document}